%% file: tokyo2023.tex
\documentclass[11pt]{amsart}

\usepackage[colorlinks = true,
            linkcolor = red,
            urlcolor  = magenta,
            citecolor = black,
            anchorcolor = blue]{hyperref}
\usepackage{color}
\usepackage[shortalphabetic,initials,nobysame]{amsrefs}
\usepackage{upgreek}
\usepackage{tikz}
\usepackage[utf8]{inputenc}
\usetikzlibrary{arrows}
\usepackage{xcolor}

\usepackage{graphicx}
\usepackage{wrapfig}

\usepackage{amssymb}

\usepackage{mathrsfs}
\usepackage[mathcal]{eucal}
\usepackage{bm}

\renewcommand{\eprint}[1]{\href{https://arxiv.org/abs/#1}{arXiv:#1}}
\newcommand{\linkedbibtitle}[1]{\IfEmptyBibField{url}{\textit{#1}}%
  {\begingroup\hypersetup{urlcolor=black}\href{\BibField{url}}{\textit{#1}}\endgroup}}
\BibSpec{misc}{%
    +{}  {\PrintPrimary} {transition}
    +{,} { \linkedbibtitle} {title}
    +{,} { \PrintDateB} {date}
    +{,} { \url} {url}
    +{.} { } {note}
}
\BibSpec{book}{%
    +{}  {\PrintPrimary} {transition}
    +{,} { \linkedbibtitle} {title}
    +{,} { } {series}
    +{,} { \voltext} {volume}
    +{,} { } {publisher}
    +{,} { } {address}
    +{,} { \PrintDateB} {date}
    +{.} { } {note}
}
\BibSpec{incollection}{%
    +{}  {\PrintAuthors} {author}
    +{,} { \linkedbibtitle} {title}
    +{,} { } {booktitle}
    +{,} { \voltext} {volume}
    +{,} { } {publisher}
    +{,} { \PrintDateB} {date}
    +{,} { pp.~} {pages}
    +{.} { } {note}
}

\newcommand{\hcL}{\hat{\cL}}

\BibSpec{article}{%
    +{}  {\PrintAuthors}                {author}
    +{,} { \linkedbibtitle}                     {title}
     +{,} {}                            {note}
    +{.} { }                            {part}
    +{:} { \textit}                     {subtitle}
    +{,} { \PrintContributions}         {contribution}
    +{.} { \PrintPartials}              {partial}
    +{,} { }                            {journal}
    +{}  { \textbf}                     {volume}
    +{}  { \PrintDate}                {date}
    +{,} { \issuetext}                  {number}
    +{,} { \eprintpages}                {pages}
    +{,} { }                            {status}
    +{,} { \DOI}                   {doi}
    +{,} { \eprint}        {eprint}
      +{,} {\publisher}                {publisher}
    +{,} { \address}                   {address}
    +{}  { \parenthesize}               {language}
    +{}  { \PrintTranslation}           {translation}
    +{;} { \PrintReprint}               {reprint}
    +{.} {}                             {transition}
    +{}  {\SentenceSpace \PrintReviews} {review}
}

\newtheorem{Thm}{Theorem}[section]

\newtheorem{Prop}[Thm]{Proposition}
\newtheorem{Cor}[Thm]{Corollary}
\newtheorem{Con}[Thm]{Conjecture}

\theoremstyle{definition}
\newtheorem{Def}[Thm]{Definition}
\theoremstyle{remark}

\newtheoremstyle{named}{}{}{\itshape}{}{\bfseries}{.}{.5em}{#1 #3}
\theoremstyle{named}

\def\g{\mathfrak{g}}
\def\Frenkel:2013uda{\mathfrak{h}}

\def\cL{\mathcal{L}}

\def\bo{\textbf{o}}

\def\=>{\Longrightarrow}

\def\to{\longrightarrow}

\def\o+{\oplus}
\def\bo+{\bigoplus}

\def\<{\langle}
\def\>{\rangle}
\def\({\left(}
\def\){\right)}

\def\^{\wedge}
\def\+{\dagger}

\def\dd[#1,#2]{\frac{d#1}{d#2}}
\def\del[#1,#2]{\frac{\partial #1}{\partial #2}}
\def\over[#1]{\overline{#1}}
\def\vec[#1]{\overrightarrow{#1}}

\makeatletter
\def\mr@ignsp#1 {\ifx\:#1\@empty\else #1\expandafter\mr@ignsp\fi}%
\newcommand{\multiref}[1]{\begingroup
\xdef\mr@no@sparg{\expandafter\mr@ignsp#1 \: }%
\def\mr@comma{}%
\@for\mr@refs:=\mr@no@sparg\do{\mr@comma\def\mr@comma{,}\ref{\mr@refs}}%
\endgroup}
\makeatother

\newcommand{\hypref}[2]{\ifx\href\asklFrenkel:2013udaas #2\else\href{#1}{#2}\fi}

\newcommand{\figref}[1]{Fig.~\multiref{#1}}

\tikzset{->-/.style={decoration={
  markings,
  mark=at position .5 with {\arrow{latex}}},postaction={decorate}}}
\tikzset{
    >=latex
    }

\newcommand{\nc}{\newcommand}
\nc{\on}{\operatorname}
\nc{\la}{\lambda}
\nc{\wh}{\widehat}
\nc{\ghat}{\wh\g}
\nc{\mb}{\mathbf}

\hypersetup{pdftitle={Quantum Geometry, Integrability, and Opers},
  pdfauthor={Peter Koroteev}}
\begin{document}
\title{Quantum Geometry, Integrability, and Opers}

\author[P. Koroteev]{Peter Koroteev}
\address{
Department of Mathematics,
University of California,
Berkeley, CA 94720, USA
and 
Department of Mathematics,
State University of New York, 
Buffalo, NY, 14068, USA
}

\date{December 29, 2023. Revised November 30, 2026}

\numberwithin{equation}{section}

\begin{abstract}
This review discusses the relations between Calogero--Ruijsenaars many-body systems, quantum spin chains, spaces of opers, double affine Hecke algebras, and enumerative invariants of quiver varieties. We describe how Wronskians of opers give Bethe equations and classical spectral levels. For the elliptic and double elliptic systems we discuss recent constructions and the questions that remain open.
{\it Contribution to the proceedings of the 16th Mathematical Society of Japan Seasonal Institute, \emph{Elliptic Integrable Systems, Representation Theory and Hypergeometric Functions}. The workshop took place at the University of Tokyo, July 31--August 4, 2023.}
\end{abstract}

\maketitle


\section{Introduction}
In this work, we review the current status of integrable models of the Calogero-Ruijsenaars family in the context of recent developments in algebraic geometry and representation theory. This study is strongly motivated by progress in string theory and supersymmetric quantum field theories beginning in the early 90s \cites{Seiberg:1994rs,Seiberg:1994aj}. 

Among supersymmetric theories with eight supercharges, the four-dimensional $\mathcal N=2^*$ theory plays a distinguished role. It is obtained from $\mathcal N=4$ Yang--Mills theory by giving the adjoint hypermultiplet a mass. Once the gauge group is fixed, the parameters of the theory include the complex gauge coupling and this mass \cite{Donagi:1995cf}.

This relative simplicity of the parameter space of the model quickly attracted experts in representation theory, algebraic geometry, and integrable systems. Physics results have been expeditiously interpreted mathematically, fostering a rich and fruitful synergy between these disciplines over the past couple of decades.

\subsection*{The Diamond}
The nine ovals in Figure~\ref{fig:diamond2023} classify the type-$A$ Calogero--Ruijsenaars systems by the dependence of their Hamiltonians on momentum and position. These variables can be elliptic, multiplicative, or additive. Starting at the double elliptic system (DELL), a step down and to the left specializes the position dependence, while a step down and to the right specializes the momentum dependence. Every path from DELL to rational Calogero--Moser specializes each variable twice. Reflection about the vertical axis exchanges the two variables.

\begin{figure}[t]
\centering
\input{oval-key}
\caption{How to read an oval in the diamond. The upper rows describe the many-body system and its oper and Bethe data in the magnetic frame. The lower rows give the algebras at finite and large particle number.}
\label{fig:notations}
\end{figure}

Figure~\ref{fig:notations} gives the order of the entries. The first row names the many-body system and gives the coordinate types for one particle, with momentum first and position second. We write $\mathcal E$ for an elliptic curve, $\mathbb C^\times$ for a multiplicative coordinate, and $\mathbb C$ for an additive coordinate. The labels describe the dependence of the Hamiltonians on a regular chart of the phase space. For $n$ particles this space has dimension $2n$ before center reduction. In an algebraically completely integrable system the $n$ independent Poisson-commuting Hamiltonians give a Lagrangian fibration on the regular locus, with fibers that compactify to abelian varieties.

We use CM for Calogero--Moser and RS for Ruijsenaars--Schneider, with the prefixes $r$, $t$, and $e$ for rational, trigonometric, and elliptic dependence. For example, the blue product $\mathbb C_p^\times\times\mathcal E_x$ beside eRS in oval~3 records that its Hamiltonians are Laurent polynomials in the multiplicative momenta and elliptic in position differences.

The next two rows give the oper and the magnetic Bethe model with its Wronskian relation. Consider tRS at the center. Here the oper is a $(G,q)$-oper with $G=SL_N$. Its proper Wronskian minors satisfy the $QQ$-system and give the XXZ Bethe equations, while its full Wronskian fixes the tRS characteristic polynomial on the canonical degree-one chart \cite{Koroteev:2022vg}. The frame matters when reading the neighboring ovals. In the magnetic frame tCM pairs with XXX and rRS with trigonometric Gaudin. In the electric frame the two spin systems exchange places. The elliptic difference and differential opers in the upper part of the diagram carry additional period data \cites{GKSZ:2024Toappear,KZEllipticPrep}.

Below the horizontal rule are the algebras at finite and large particle number. In oval~5 the finite-rank entry is the double affine Hecke algebra (DAHA). Its spherical subalgebra quantizes the multiplicative Calogero--Moser space at fixed puncture conjugacy class \cites{Oblomkov:aa,Cherednik-book}. At large rank the stable spherical DAHA has an elliptic Hall algebra realization. With the appropriate central parameters, this connects the quantum toroidal algebra of $\mathfrak{gl}_1$ to equivariant K-theory of Hilbert schemes of points on $\mathbb C^2$ \cites{Matsuo:2023aa,Koroteev:2020mxs}. In the figure, tCA and rCA denote trigonometric and rational Cherednik algebras, DIM denotes Ding--Iohara--Miki, and DDCA denotes a deformed double current algebra. A question mark in either algebra row records an open identification.

\subsection*{Enumerative Geometry.}
The two labeled diagonals in Figure~\ref{fig:diamond2023} organize the Laumon-space constructions. The ordinary Laumon diagonal runs through ovals $2\to5\to8$. At its last two vertices the tRS and tCM wave functions admit descriptions in equivariant K-theory and cohomology of ordinary Laumon spaces, respectively. The affine Laumon diagonal runs through $1\to3\to6$, from DELL through eRS to eCM. The stationary limit and its normalization enter the elliptic eigenfunction problem. Section~\ref{Sec:QuantumManyBody} gives the status of the enumerative identification at each of these three vertices.

\subsection*{Supersymmetric Gauge Theories.} 
The relation of Seiberg--Witten geometry to integrable systems was developed in \cites{Gorsky:1995zq,Gorsky:2000px,Gorsky:1997mw}. In four, five, and six dimensions, theories with adjoint matter and a codimension-two defect lead to wave functions for eCM, eRS, and DELL, respectively. The five-dimensional theory is compactified on a circle and its six-dimensional counterpart on a torus. The precise enumerative statement depends on the dimension. In particular, the DELL eigenfunction identity for the elliptic genus is a conjecture \cite{Koroteev:2019gqi}.

\subsection*{Reading the elliptic ovals.}
On the edge $3\to6$ one fixes the spectral eigenline of the particle matrix and an ordered divisor at which to evaluate its sections. At oval~3 one takes shifted evaluations of a Baker--Akhiezer section. In the small-step limit these become its jets at oval~6, where a one-point differential oper gives another description. The accompanying finite Bethe models act on the same space of color assignments at fixed sites, with each color occurring once. The first-order terms of the normalized inverse Felder transports give the Kronecker permutation operators. At oval~2 the XYZ entry refers to Bethe equations with a choice of boundary character. Section~\ref{Sec:MiddleOval} explains the extra data needed to realize a prescribed vertex chain. At oval~1 the cyclic theta construction gives a surface $QQ$-system and Bethe equations. A microscopic lattice Hamiltonian with these equations remains to be constructed.

We begin in Section~\ref{Sec:MiddleOval} with the central oval, where both the Bethe equations and the classical spectral level can be read directly from the Wronskians. Section~\ref{Sec:QuantumManyBody} takes up the quantization of the particle systems and their enumerative interpretation.

\subsection*{Acknowledgments} 
I am grateful to A. Zeitlin, E. Frenkel, and D. Sage for their continuing collaboration and encouragement. I maintain an \href{https://math.berkeley.edu/~pkoroteev/diamond.html}{interactive version of the diamond} with references for individual ovals and arrows, and welcome additions from colleagues. I thank H. Konno and J. Shiraishi for inviting me to present this work in Tokyo in August 2023.

\clearpage
\begin{figure}[p]
\centering
\makebox[\linewidth][c]{\includegraphics[width=1.075\linewidth]{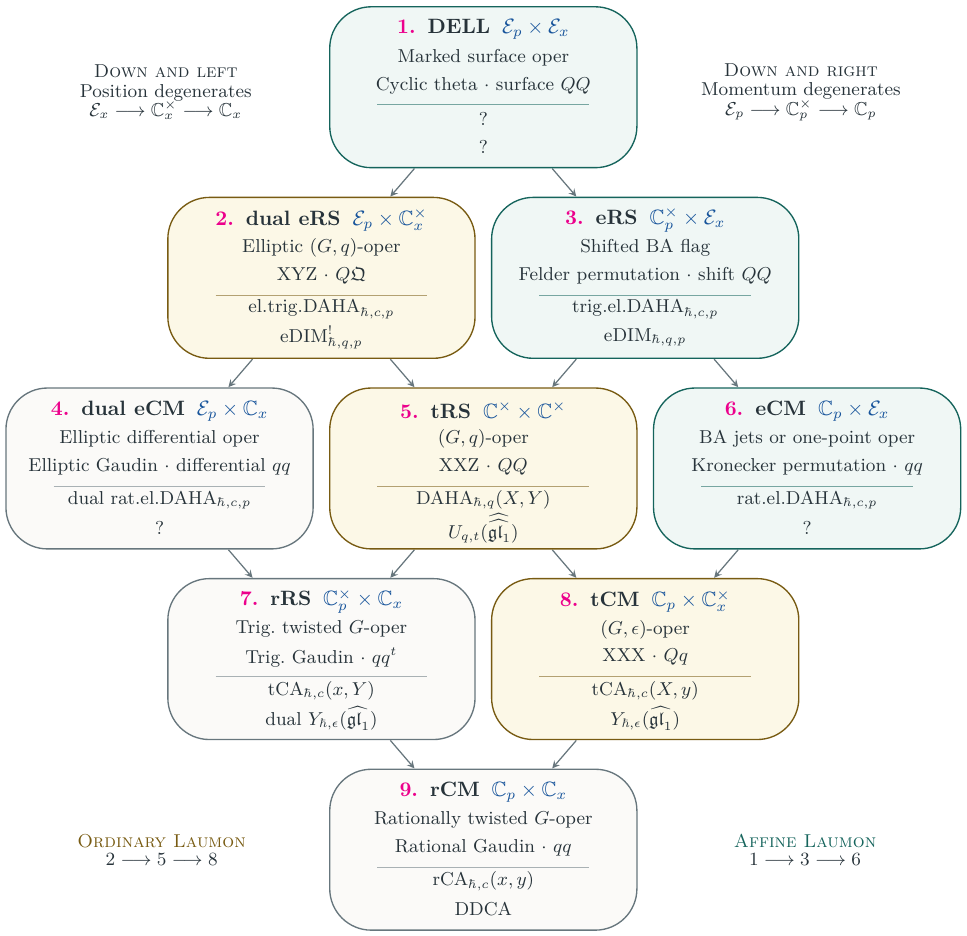}}
\caption{The type-$A$ diamond of integrability in the magnetic frame. The arrows specialize one coordinate at a time. At oval~2 the elliptic algebra entries are proposed and $!$ denotes the mirror polarization. The dual elliptic algebra at oval~4 is also proposed. At oval~7, \emph{dual} denotes the rational difference realization of the affine Yangian. In the algebra subscripts $p$ is the elliptic parameter. Some of the new results announced here will appear shortly in \cite{DiamondPrep}.}
\label{fig:diamond2023}
\end{figure}
\clearpage

\section{\texorpdfstring{$(G,q)$}{(G,q)}-Opers, tRS model, and XXZ Chain}\label{Sec:MiddleOval}
The space of opers provides a framework for relating integrable systems and their dualities \cites{Frenkel:2020,KZMinors,Koroteev:2022vg}. Here we review type $A$ with periodic twist. With M. Shim and R. Singh we studied reflection-invariant $q$-opers for open chains and the folding of type-$A$ data to the classical-group Ruijsenaars systems \cites{KSSOpen,KSSBeyond}.

We begin with the middle oval of the diamond. Throughout this section set $N=r+1$, let $X=\mathbb{P}^1\backslash\{0,\infty\}$, and let $M_q(z)=qz$, where $q\in\mathbb C^\times$ is not a root of unity. A superscript $q$ denotes pullback by $M_q$. We use the inverse-transport convention for a $q$-connection: $A(z):E_{qz}\to E_z$. Its inverse is the forward-transport convention $E_z\to E_{qz}$; these describe the same geometric data. This choice makes the shifted Wronskians below literal determinants of vectors in one fiber.


\begin{Def}    \label{qopflagD}
Let $U \subset X$ be a Zariski open dense subset and let $V=U \cap M_q^{-1}(U)$.
 A meromorphic $(GL(r+1),q)$-{\em oper} on $\mathbb{P}^1$ is a triple $(A, E, \mathcal{L}_{\bullet})$, where $E$ is a vector bundle of rank $r+1$ on $\mathbb{P}^1$ and $\mathcal{L}_{\bullet}$ is the corresponding complete flag of the vector bundles, 
  $$\mathcal{L}_{r+1}\subset ...\subset \mathcal{L}_{i+1}\subset\mathcal{L}_i\subset\mathcal{L}_{i-1}\subset...\subset \mathcal{L}_1=E,$$ 
  where $\mathcal{L}_{r+1}$ is a line bundle and $\mathcal L_{r+2}=0$, so that the meromorphic $(GL(r+1),q)$-connection
  $A:E^q\to E$
  satisfies the following conditions:\\ 
i) $A\cdot\mathcal L_i^q\subset\mathcal L_{i-1}$ for $i=2,\dots,r+1$,\\
ii) $A$ is invertible on $V$ and the induced maps
  $$\bar A_i:\mathcal L_i^q/\mathcal L_{i+1}^q\to\mathcal L_{i-1}/\mathcal L_i,\qquad i=2,\dots,r+1$$ are isomorphisms on $V$. \\
  An $(SL(r+1),q)$-oper carries, in addition, a trivialization of $\det E$ in which $\det A=1$ on $V$.
\end{Def}

A change of trivialization by $g(z)\in SL(r+1)(\mathbb C(z))$ acts on the connection matrix by
\begin{equation}   
\label{gauge tr}
A(z)\mapsto g(z)A(z)g(qz)^{-1}.
\end{equation}

The above definition can be reformulated in a local form -- given a section $s(z)$ of $\mathcal{L}_{r+1}$, the oper condition is as follows. Consider the following determinants
\begin{equation}\label{altqW} 
W_k(s)(z)=s(z)\wedge A(z)s(qz)\wedge\cdots\wedge
 A(z)A(qz)\cdots A(q^{k-2}z)s(q^{k-1}z),\quad k=2,\dots,r+1.
\end{equation}
These are sections of $\bigwedge^k\mathcal L_{r-k+2}$; after trivializing this line bundle we regard them as scalar functions. The oper conditions require them to be nonzero on the open set under consideration.
\begin{Def} 
We say that the oper has regular singularities specified by nonzero polynomials $\Lambda_1(z),\dots,\Lambda_r(z)$ when, in compatible local trivializations, the zero divisor of $\bar A_{r+2-i}$ is the divisor of $\Lambda_i$, for $i=1,\dots,r$. Multiplicities are included. We number these data starting from the line $\mathcal L_{r+1}$; this is also the ordering of the minors and $QQ$-relations below.
\end{Def}
In local terms, the above definition implies the following
\begin{align}
\label{eq:WPDefs}
F_0=F_1=1,\qquad
F_k(z)=\prod_{a=1}^{k-1}\prod_{b=0}^{k-a-1}\Lambda_a(q^b z),\qquad k=2,\dots,r+1.
\end{align}
Here $F_k$ records the prescribed singularity factor of the $k$-th Wronskian, up to a nonzero scalar normalization. In particular,
$F_{i-1}(qz)F_{i+1}(z)=\Lambda_i(z)F_i(z)F_i(qz)$.

\begin{Def}
The $(SL(r+1),q)$-oper is called $Z-twisted$ if there exists $g(z)\in SL(r+1)(z)$ such that
\begin{eqnarray}    \label{Ag}
A(z)=g(z)Z g(qz)^{-1},
\end{eqnarray}
where $Z$ is a diagonal element of $SL(r+1)$.
\end{Def}

\begin{Def}  
A Miura $(SL(r+1),q)$-oper is a quadruple $(A,E, \mathcal{L}_{\bullet}, \hcL_\bullet)$, where 
$(A,E, \mathcal{L}_{\bullet})$ triple is $(SL(r+1),q)$-oper and the complete flag $\hcL_\bullet$ of subbundles in $E$ is preserved by the q-connection $A$.  
\end{Def}

We work with generic twist parameters: $Z=\operatorname{diag}(\xi_1,\dots,\xi_N)$ with $\xi_i/\xi_j\notin q^{\mathbb Z}$ for $i\ne j$. In particular, $Z$ is regular semisimple, meaning that its eigenvalues are distinct. In this nonresonant setting the $N!$ orderings of its eigenlines give the $N!$ invariant complete flags. Fix one ordering for the discussion below.

Relative position between flags $\mathcal{L}_\bullet$ and $\hat{\mathcal{L}}_\bullet$  can be described using the following determinant
\begin{equation}\label{qD}
\mathcal D_k(s)=s(z)\wedge Zs(qz)\wedge\cdots\wedge Z^{k-1}s(q^{k-1}z)
\wedge e_{k+1}\wedge\cdots\wedge e_N,
\end{equation}
for $k=1,\dots,N$, with $\mathcal D_0=1$. Here we have chosen the invariant flag so that these are the leading minors; $s(z)$ is a polynomial representative of a local generator of $\mathcal L_N$, with components $s_1(z),\dots,s_N(z)$ in the twist eigenbasis. Changing from leading to trailing minors reverses the node ordering, and must also reverse the twist and singularity indices.
The functions $\mathcal{D}_k(s)$ have a subset of zeros, which coincide with those of
$F_k(z)$. The remaining zeros of $\mathcal{D}_k(s)$ are given by points at which the two flags fail to be in a generic position. We have accordingly
\begin{equation}
  \underset{1\leq i,j\leq k}{\det}\left[\xi_i^{j-1}s_i(q^{j-1}z)\right]
  =\beta_k F_k(z)\mathcal V_k(z),
\label{eq:MiuraDetForm}
\end{equation}
where $i,j = 1,\dots,k$ and $s_l^{(m)}=s_l(q^mz)$,  $Z={\rm diag}(\xi_1, \dots, \xi_{r+1})$, and 
\begin{equation} 
\mathcal V_k(z)=\prod_{a=1}^{m_k}(z-s_{k,a}),\qquad m_k=\deg\mathcal V_k,
\label{eq:BaxterRho}
\end{equation}
where $\beta_k$ are nonzero constants. The boundary polynomials are $\mathcal V_0=\mathcal V_N=1$: the full determinant has no zeros beyond its prescribed singularities. The complex numbers $s_{k,a}$ are the Bethe roots, not the components $s_i(z)$ of the oper section.

Before stating the correspondence, fix $q$, the ordered twist $Z$, the singularity polynomials $\Lambda_i$, and the degrees $m_i$ for $i=1,\dots,r$. Set $Q_i^+=\mathcal V_i$, with $Q_0^+=Q_N^+=1$. A polynomial solution of the $QQ$-system means a tuple of monic polynomials $Q_i^+$ of these degrees and auxiliary polynomials $Q_i^-$ satisfying the displayed identities coefficient by coefficient in $z$. Thus the unknowns are polynomial coefficients, or equivalently the unordered roots $s_{i,k}$; the spectral variable $z$ and the parameters $q,\xi_i,\Lambda_i$ are not unknowns to be solved for.

We restrict to the nondegenerate locus. Here the roots of each $Q_i^+$ are simple and nonzero, no two roots of the same or adjacent $Q^+$-polynomials lie in the same $q^{\mathbb Z}$-orbit, and the roots of $Q_i^+$ are $q$-distinct from the zeros of $\Lambda_i$. Two nonzero numbers $u,v$ are called $q$-distinct if $u/v\notin q^{\mathbb Z}$. These genericity conditions, together with the nonresonant twist above, in particular make all divisions in the Bethe equations below legitimate.

The remarkable property of polynomials \eqref{eq:BaxterRho} is that they solve the $QQ$-system.
\begin{Thm}[\cite{KSZ}]\label{qWthKSZ}
Polynomials $\{\mathcal{V}_k(z)\}_{k=1,\dots, r}$ provide the solution to the $QQ$-system 
\begin{equation}
 \xi_i Q_i^+(qz)Q_i^-(z)-\xi_{i+1}Q_i^+(z)Q_i^-(qz)
 =(\xi_i-\xi_{i+1})\Lambda_i(z)Q_{i-1}^+(qz)Q_{i+1}^+(z),
\label{eq:dQQrelations}
\end{equation}
so that $Q^+_j(z)=\mathcal{V}_j(z)$. The polynomials $Q^+_j, Q^-_j$ for $j=1,\dots,r$ can be presented using minors
\begin{equation}
Q^+_j(z)\ \propto\ \frac{1}{F_j(z)}\frac{\det M_{1,\ldots,j}(z)}{\det V_{1,\ldots,j}},
\qquad
Q^-_j(z)\ \propto\ \frac{1}{F_j(z)}\frac{\det M_{1,\ldots,j-1,j+1}(z)}{\det V_{1,\ldots,j-1,j+1}},
\label{eq:QPolyM}
\end{equation}
where $F_j$ is defined in \eqref{eq:WPDefs}, and
\begin{equation}
M_{i_1,\ldots,i_j} = \begin{bmatrix} \,  s_{i_1} & \xi_{i_1} s_{i_1}^{(1)} & \cdots & \xi_{i_1}^{j-1} s_{i_1}^{(j-1)} \\ \vdots & \vdots & \ddots & \vdots \\  s_{i_j}  & \xi_{i_j} s_{i_j}^{(1)} & \cdots & \xi_{i_j}^{j-1} s_{i_j}^{(j-1)} \, \end{bmatrix}\,,
\qquad
V_{i_1,\ldots,i_j} =\begin{bmatrix} \,  1 & \xi_{i_1}  & \cdots & \xi_{i_1}^{j-1} \\ \vdots & \vdots & \ddots & \vdots \\  1 & \xi_{i_j} & \cdots &\xi_{i_j}^{j-1}\, \end{bmatrix}\,,
\label{eq:MM0Ind2}
\end{equation}
are the quantum Wronskian and the Vandermonde matrix respectively. The symbol $\propto$ means equality up to a nonzero constant, not a function of $z$. Normalize $Q_j^+$ to be monic and choose the constant in $Q_j^-$ so that \eqref{eq:dQQrelations} holds. The determinant identity underlying this formula is
\[
\xi_j d_j(qz)\widetilde d_j(z)-\xi_{j+1}d_j(z)\widetilde d_j(qz)
=(\xi_j-\xi_{j+1})d_{j-1}(qz)d_{j+1}(z),
\]
where $d_j=\det M_{1,\ldots,j}/\det V_{1,\ldots,j}$, $\widetilde d_j$ uses the row set $1,\ldots,j-1,j+1$, and $d_0=1$. Dividing by the singularity factors gives \eqref{eq:dQQrelations}.
\end{Thm}

We can now combine the correspondences of \cites{KSZ,Frenkel:2020,Koroteev:2020mxs,KZMinors}.
\begin{Thm}\label{Th:3wayXXZ}
For the fixed data and nondegeneracy conditions specified above, there are natural bijections between the following sets:

i) Polynomial solutions $(Q_i^+,Q_i^-)_{i=1}^r$ of \eqref{eq:dQQrelations}, with $\deg Q_i^+=m_i$ and $Q_i^+$ monic.

ii) Nondegenerate solutions in $s_{i,k}$, modulo permutations of the roots at each fixed node $i$, of the Bethe Ansatz equations for the $\mathfrak{sl}(r+1)$ XXZ spin chain:
\begin{equation}    \label{eq:bethe}
\frac{Q^+_{i}(qs_{i,k})}{Q^+_{i}(q^{-1}s_{i,k})} \frac{\xi_i}{\xi_{i+1}}=
-\frac{\Lambda_i(s_{i,k})Q^+_{i-1}(qs_{i,k})Q^+_{i+1}(s_{i,k})}
{\Lambda_i(q^{-1}s_{i,k})Q^+_{i-1}(s_{i,k})Q^+_{i+1}(q^{-1}s_{i,k})},
\end{equation}
where $i=1,\ldots,r$ and $k=1,\ldots,m_i$.

iii) Isomorphism classes of nondegenerate $Z$-twisted Miura $(SL(r+1),q)$-opers with these singularity data and these degrees of the residual minors.
\end{Thm}
The passage from i) to ii) evaluates \eqref{eq:dQQrelations} at $z=s_{i,k}$ and $z=q^{-1}s_{i,k}$ and eliminates $Q_i^-(s_{i,k})$. Conversely, the nondegenerate Bethe equations give the polynomial $Q_i^-$ under the nonresonance assumption \cites{KSZ,Frenkel:2020}. A space of opers and its function algebra are different objects; below $\operatorname{Fun}$ denotes the latter. When using roots as algebra generators, we take functions invariant under their permutations, or use the coefficients of the monic $Q_i^+$ instead.

%

\subsection{The tRS model and \texorpdfstring{$(SL(N),q)$}{(SL(N),q)}-opers}
We now restrict to the \textit{canonical degree-one locus}; this is an additional restriction, not a property of every $q$-oper. Choose a polynomial representative of the line $\mathcal L_N$ whose components all have degree one, with distinct twist eigenvalues, and require that the only prescribed Wronskian zeros are those of the full determinant \cites{KSZ,Koroteev:2022vg}. In the leading-minor convention this means $\Lambda_i=1$ for $i<r$ and $\Lambda_r=\Lambda$, where $\Lambda(z)=\prod_{i=1}^N(z-a_i)$ has distinct nonzero roots. Equivalently, on this generic degree-one chart the normalized $k$-th minors have degree $k$, and the normalized full minor equals $\Lambda$. The corresponding Bethe degrees are $m_i=i$ for $i=1,\dots,r$.

Nonzero leading coefficients of the components can be removed by diagonal changes of the twist eigenbasis and by rescaling the generator of the line. Hence we may write $s_i(z)=z-p_i$. Conversely, prescribing this degree-one form and the full determinant equation gives precisely the chart used here. We denote this canonical oper space by $q\operatorname{Op}_Z^\Lambda$ below. For the $N$-particle formulas we retain the scalar $GL(N)$ extension of the $SL(N)$ data, so that both ordered parameter sets $\xi_i$ and $a_i$ can be kept independent; determinant and center-of-mass constraints may subsequently be imposed for the $SL(N)$ reduction.

For clarity, let $K=\mathbb C(q,\xi_1,\dots,\xi_N,a_1,\dots,a_N)$ be the field of generic \emph{parameters}, and let $R_p=K[p_1^{\pm1},\dots,p_N^{\pm1}]$. Let $R_s$ be the $K$-algebra of polynomials in the Bethe roots, invariant under $\prod_i S_{m_i}$, localized at the denominators used in the Bethe equations and at the root discriminants. All subsequent coordinate-ring identities are understood on the indicated generic charts. In particular, we never quotient a field containing the unknown momenta by a nonzero relation. Theorem~\ref{Th:3wayXXZ} gives
\begin{eqnarray}
\operatorname{Fun}(q\operatorname{Op}_Z^\Lambda)
\cong R_p/I_{\mathrm{Wr}}\cong R_s/I_{\mathrm{Bethe}},
\end{eqnarray}
where $I_{\mathrm{Wr}}$ equates the coefficients of the normalized full determinant and $\Lambda(z)$. The ideal $I_{\mathrm{Bethe}}$ is the intersection with $R_s$ of the permutation-stable ideal generated by the denominator-cleared equations \eqref{eq:bethe} in the localized root algebra, in this degree sector.

The variables $p_i\in\mathbb C^\times$ are \emph{multiplicative}, or exponentiated, tRS momenta. Locally, $x_i=e^{X_i}$ and $p_i=e^{P_i}$ give the additive Darboux variables $X_i,P_i$, with symplectic form $\sum_i dP_i\wedge dX_i$. They are not additional Bethe roots: the latter are the roots of the larger minors $Q_i^+$. In our degree-one normalization, set $\widehat V_{ij}=(q\xi_i)^{j-1}$ for $1\leq i,j\leq N$. Then
\begin{equation}
M_{1,\ldots,N}(z)=z\widehat V+M_{1,\ldots,N}(0).
\end{equation}
The normalized full Wronskian is
\begin{equation}
P(z)=\frac{\det M_{1,\ldots,N}(z)}{\det\widehat V}=\Lambda(z),
\end{equation}
or, after inverting $\widehat V$, the formula reads
\begin{equation}
P(z)=\text{det}\left(z-T\right)\,, 
\label{eq:BaxterFlavorPol}
\end{equation}
where 
\begin{equation}
T=-M_{1,\ldots,N}(0)\widehat V^{-1}.
\label{eq:tRSLaxMatrix}
\end{equation}
As it was shown in \cite{KSZ} the latter is the Lax matrix for the tRS model
\begin{eqnarray}
T_{ij}=  \frac{\prod\limits_{m \neq j}^{r+1} \left(q^{-1} \xi_i \,  -  \xi_m \,  \right) }{\prod\limits_{l\neq j}^{r+1}(\xi_j-\xi_l)}   p_i\,, \qquad i,j=1,\dots,r+1\,
\end{eqnarray}
whose characteristic polynomial \eqref{eq:BaxterFlavorPol} generates integrals of motion
\begin{equation}\label{eq:tRSW3}
\det(zI_N-T)=\sum_{k=0}^N(-1)^k H_k^{tRS}(q,\xi,p)z^{N-k},\qquad H_0^{tRS}=1,
\end{equation}
such that
\begin{equation}\label{eq:tRSmagdef}
H^{tRS}_k=q^{-k(k-1)/2}\sum_{\substack{\mathcal I\subset\{1,\dots,N\}\\|\mathcal I|=k}}
\prod_{\substack{i\in\mathcal I\\j\notin\mathcal I}}
\frac{q^{-1}\xi_i-\xi_j}{\xi_i-\xi_j}\prod_{m\in\mathcal I}p_m.
\end{equation}
The factor $q^{-k(k-1)/2}$ follows from the principal minors of the displayed Lax matrix and fixes our Hamiltonian normalization. For example, $\det T=q^{-N(N-1)/2}\prod_i p_i$. Equating the characteristic polynomial to $\Lambda$ gives
\begin{Thm}[\cite{KSZ}]\label{Th:XXz/tRS}
There is an isomorphism of algebras 
\begin{equation}
{\rm Fun}\big(q{\rm Op}_Z^{\Lambda}\big)\cong \frac{R_p}
{\left( H^{tRS}_k-e_k(a_1, \dots, a_{r+1})\right)_{k=1,\dots, r+1}},
\end{equation}
where $e_k(a_1,\dots,a_N)=\sum_{|I|=k}\prod_{i\in I}a_i$ is the elementary symmetric polynomial of degree $k$. Here $q,\xi_i,a_i$ are fixed and the equations determine the momenta $p_i$.
\end{Thm}

\subsection{Quantum/Classical Duality -- Magnetic Frame}
Here \textit{quantum/classical duality} identifies the Bethe eigenvalue problem of the quantum XXZ chain with a spectral-level-set problem for a \emph{classical} tRS system; it does not yet quantize the tRS phase space. The geometric formulation is developed in \cite{Gaiotto_2013}; for the corresponding XXX/rRS and Gaudin/rCM constructions see \cite{GZZSpin}. On a tRS Darboux chart the symplectic form is $\Omega=\sum_{i=1}^N d\log p_i\wedge d\log x_i$, and the Hamiltonians Poisson commute. A \textit{frame} specifies which fixed spin-chain parameters become classical particle positions and which become the prescribed Lax eigenvalues. The magnetic frame uses $x_i=\xi_i$ and $\operatorname{Spec}(T)=\{a_i\}$; the electric frame below interchanges these roles. The momenta in the two frames are generally different functions on the same correspondence space.

\begin{Cor}\label{tqqc}
We have the following isomorphisms:
\begin{eqnarray}
R_s/I_{\mathrm{Bethe}}\cong {\rm Fun}\big(q{\rm Op}_Z^{\Lambda}\big)\cong \frac{R_p}
{\left(\{H^{tRS}_k-e_k(\{a_i\})\}_{k=1,\dots, r+1}\right)}\,.
\end{eqnarray}
The latter algebra is the coordinate algebra of the scheme-theoretic intersection $\mathscr L_1\cap\mathscr L_2$ of two Lagrangian subvarieties on the regular locus:
\begin{eqnarray} 
\mathscr{L}_1:\,\,\{x_i=\xi_i\}, \quad \mathscr{L}_2:\,\,\{H^{tRS}_i=e_i(\{a_i\})\}\,,\quad i=1,\dots, r+1\,.
\end{eqnarray}
\end{Cor}

The above theorem establishes the quantum/classical duality in the \textit{magnetic frame} according to the terminology of \cite{Koroteev:2023aa} since the tRS coordinates $x_i$ are represented in terms of twist parameters $\xi_i$ of the $q$-oper.

\subsection{Quantum/Classical Duality -- Electric Frame}
In the \textit{electric frame}, the positions are the singularity parameters $x_i=a_i$, while the Lax eigenvalues are the twist parameters $\xi_i$ (see Section 4 of \cite{Koroteev:2023aa}). Write $p'_i$ for its multiplicative momenta to distinguish them from the magnetic $p_i$, and put $R_{p'}=K[(p'_1)^{\pm1},\dots,(p'_N)^{\pm1}]$. The interaction parameter is inverted. With the same characteristic-polynomial normalization, the spectral equations are
\begin{equation}\label{eq:tRSEqnsElectric}
H^{!,tRS}_k:=q^{k(k-1)/2}\sum_{\substack{\mathcal I\subset\{1,\dots,N\}\\|\mathcal I|=k}}
\prod_{\substack{i\in\mathcal I\\j\notin\mathcal I}}
\frac{qa_i-a_j}{a_i-a_j}\prod_{m\in\mathcal I}p'_m
=e_k(\xi_1,\dots,\xi_N),\quad k=1,\dots,N.
\end{equation}

\begin{Thm}[\cite{Koroteev:2023aa}]\label{Th:tRSXXZEl}
There is an isomorphism of algebras:
\begin{eqnarray}\label{eq:FunAlgFl}
{\rm Fun}\big(q{\rm Op}_Z^{\Lambda}\big)\cong
\frac{R_{p'}}{(H_k^{!,tRS}-e_k(\xi_1,\dots,\xi_N))_{k=1}^N}.
\end{eqnarray}
The latter algebra is the coordinate algebra of the corresponding intersection in the electric chart, whose symplectic form is $\Omega'=\sum_i d\log p'_i\wedge d\log x_i$:
\begin{eqnarray} 
\mathscr{L}_1:\,\,\{x_i=a_i\}, \quad \mathscr{L}_2:\,\,\{H^{!\,tRS}_i=e_i(\{\xi_i\})\}\,,\quad i=1,\dots, r+1\,,
\end{eqnarray}
where $H^{!\,tRS}_i$ are the left hand sides of \eqref{eq:tRSEqnsElectric}.
\end{Thm}

\subsection{3d Mirror Symmetry}
The exchange of particle positions with Lax spectral data is the classical action--angle duality developed by Ruijsenaars \cite{RuijAction}. The relation between the two tRS descriptions here is its complex-algebraic manifestation, with the inverse-coupling convention made explicit. We do not claim a new classical tRS self-duality. The oper and quantum K-theory descriptions identify how this known duality acts on the associated Bethe and geometric data.

In the gauge-theoretic terminology, \textit{3d mirror symmetry} exchanges Higgs- and Coulomb-branch descriptions; in the examples considered here it exchanges equivariant and K\"ahler parameters. Its manifestation in the present spectral problem exchanges the singularity and twist parameters $a_i$ and $\xi_i$ and sends $q$ to $q^{-1}$. The term \textit{bispectral duality} refers to the corresponding exchange of the two spectral problems, not to an identification of every aspect of 3d mirror symmetry with classical mechanics.

\begin{Prop}[\cite{Koroteev:2022vg}]\label{tqqcm2}
The following algebras are isomorphic
\begin{eqnarray}
\frac{R_p}{(H^{tRS}_k-e_k(a_1,\dots,a_N))_{k=1}^N}
\cong\frac{R_{p'}}{(H^{!,tRS}_k-e_k(\xi_1,\dots,\xi_N))_{k=1}^N}.
\end{eqnarray}
\end{Prop}

\begin{Thm}
The canonical degree-one oper family is self-dual under this parameter exchange, namely:
\begin{eqnarray}
{\rm Fun}\big(q^{-1}{\rm Op}_{Z'}^{\Lambda'}\big)\cong {\rm Fun}\big(q{\rm Op}_{Z}^{\Lambda}\big).
\end{eqnarray}
where $Z'=\text{diag}(a_1,\dots,a_{r+1})$ and $\Lambda'(z)=\prod_{i=1}^{r+1}(z-\xi_i)$.
\end{Thm}

See \cite{Koroteev:2023aa} for the oper and quiver-variety interpretation and further references. In the diamond \figref{fig:diamond2023} this exchange is represented by reflection about the vertical line passing through ovals 1, 5, and 9. Thus tCM is dual to rRS, XXX is dual to tGaudin, and tRS, XXZ, rCM, and rGaudin are self-dual in the indicated sense \cites{Feher:2009uw,Mukhin:2009uh,Feher:2010va,Feher:2012wc,Mukhin:wg}.
Some models in the diamond are defined as dual to known models, such as dual eCM and dual eRS. The latter two can be also obtained via degenerate limits from DELL which we shall address later.

\subsection{Limits and Degenerations}
The bottom four ovals (\#5,7,8, and 9) of the diamond have been studied in \cite{Koroteev:2022vg}. In particular, the author together with A. Zeitlin studied spaces of $(G,\epsilon)$-opers, trigonometrically twisted $G$-opers (t$G$-opers), and rationally twisted $G$-opers (r$G$-opers). For the $(G,\epsilon)$-opers, one picks $X=\mathbb{P}^1$ in Definition \ref{qopflagD} and replaces $M_q$ with the additive shift $M_\epsilon: z\mapsto z+\epsilon$. As for the t$G$-opers and r$G$-opers one considers a vector bundle $E$ together with a complete flag of subbundles $\mathcal{L}_\bullet$ over $\mathbb{P}^1$ with a meromorphic connection $\nabla: E\to E\otimes K$ satisfying $\nabla\mathcal L_i\subset\mathcal L_{i-1}\otimes K$ away from the poles, for $2\le i\le N$. One requires that the induced maps $\mathcal{L}_i/\mathcal{L}_{i+1}\to \mathcal{L}_{i-1}/\mathcal{L}_{i}\otimes K$ are isomorphisms away, perhaps, from a set of isolated singularities. The rationally twisted connection can be brought to the form $\partial_z+Z$, while the trigonometrically twisted connection can be brought to $\partial_z+
\frac{Z}{z}$. Equivalently, the latter is $\partial_\zeta+Z$ in the cylinder coordinate $z=e^\zeta$. The gauge transformation lies in $G(z)$, or in $B_+(z)$ for a twisted Miura oper; see \cites{Brinson:2021aa,Koroteev:2022vg}.

Each type of opers leads to its own set of quantum/classical dual models in the magnetic frame:
$$
\text{rRS/tGaudin}\qquad \qquad \text{tCM/XXX}\qquad \qquad\text{rCM/rGaudin}
$$
In the electric frame, the dualities go as follows
$$
\text{rRS/XXX}\qquad \qquad \text{tCM/tGaudin}\qquad \qquad\text{rCM/rGaudin}
$$
As was pointed out earlier the tRS model and the XXZ chain are self-dual under 3d mirror symmetry (bispectral duality).

The tRS, rRS, tCM, and rCM models can be obtained from the Calogero-Moser space description by taking certain limits, see \cite{Koroteev:2022vg} for review.

\subsection{Elliptic Opers}
In oval 2 the base is the elliptic curve $\mathcal E_w=\mathbb C^\times/w^{\mathbb Z}$, $0<|w|<1$, with translation $z\mapsto qz$. The oper consists of a bundle, a meromorphic difference connection and a full flag as in Definition~\ref{qopflagD}. A Miura structure supplies an invariant flag. On the split chart fix the ordered flat summands and their difference lifts. In a theta frame the components $s_i(z)$ of the cyclic line have prescribed degree and multipliers. The first $k$ translates give the Wronskian minors of \eqref{eq:QPolyM} with theta sections in place of polynomials. The $QQ$ identities follow from the same determinant identity. At separated zeros they give elliptic Bethe equations. Their divisors satisfy an additional condition: the Abel sum in each degree is fixed by the multiplier of the corresponding line bundle \cites{GKSZ:2024Toappear,KZEllipticPrep}.

When each $s_i$ has degree one, the top determinant gives the characteristic current of dual eRS in the magnetic frame, up to the fixed row and Vandermonde factors. Prescribing its section up to scale fixes the classical level on this chart. The split reconstruction uses fixed characteristics, nonresonance and separation of the shifted zeros. These hypotheses matter in comparing geometric Bethe data to a transfer matrix.

For $SL(2)$ there is a dynamical Bethe ansatz for Felder's elliptic quantum group \cite{FVel}. The transfer matrix acts on functions of a dynamical parameter in the zero-weight sector. To obtain a correspondence with an ordinary vertex XYZ chain one must also specify the representation, the Bloch multipliers and the boundary condition. The identification for a fixed vertex chain from prescribed oper data remains a question. In higher rank Hou, Sasaki and Yang give a nested Bethe ansatz for elliptic quantum-group modules and identify specific Belavin vertex-model boundary sectors \cite{HouNested}. In the mirror frame the exchange with ordinary eRS must carry the same boundary data.

\section{Enumerative Geometry and Quantum Many-Body Systems}\label{Sec:QuantumManyBody}
The next step in our construction is to quantize each classical many-body system from the diamond and find connections with enumerative geometry. We begin with the quantum equivariant K-theory of quiver varieties and the tRS system (oval \#5). Below we shall introduce quantum K-theory vertex functions and formulate difference equations that they satisfy. More details are provided in \cites{Koroteev:2018azn,Koroteev:2021aa,KZGaugeReview}.

We now use $q$ for the quantum shift of the many-body difference operators and $\hbar$ for the cotangent-scaling parameter. The classical magnetic Hamiltonians of Section~\ref{Sec:MiddleOval} are recovered with that section's oper parameter equal to $\hbar^{-1}$. The variables $z_i$ below are K\"ahler variables.

For a given quiver variety $X$, the quantum multiplication in its localized $T$-equivariant K-theory $K_T(X)[[z]]$ can be defined by studying moduli spaces of quasimaps $\textsf{QM}$ \cite{Ciocan-Fontanine:2011tg} in all positive degrees from a base curve into the variety.
\begin{equation*}
A\circledast B = A\otimes B + \sum\limits_{\textbf{d}>0}\,z^{\textbf{d}} A\circledast_{\textbf{d}} B\,.
\end{equation*}
This product enables us to construct a unital quantum K-theory ring. The quasimap moduli spaces from $\mathcal{C}$ into $X$ have a natural action of the maximal torus, lifting its action from $X$. When there are at most two special points and the base curve $\mathcal{C}$ is ${\mathbb{P}}^1$, we extend $T$ by additional torus $\mathbb{C}^{\times}_q$, which scales ${\mathbb{P}}^1$ such that the tangent space $T_{0} {\mathbb{P}}^1$ has character denoted by $q$. We will include this action in the full torus by $\mathsf{T}=T\times \mathbb{C}^{\times}_q$. We assume that the two fixed points of $\mathbb{C}^{\times}_q$ are $0$ and $\infty$.
As explained in \cite{Okounkov:2015aa}, one can construct various enumerative invariants of $X$ using virtual structure sheaves $\mathscr{O}_{\text{vir}}$ for $\textsf{QM}^{\textbf{d}}$. Using the above two marked points, one can define a vertex function (analog of Givental's J-function).
Here a quasimap is allowed finitely many base points where it fails to define a map to the GIT quotient $X$; ``nonsingular at $p_2$'' excludes a base point there. In the following formula $p_1=0$, $p_2=\infty$, $\mathrm{ev}_{p_2}$ is evaluation at $p_2$, $\mathscr V_i$ are the universal tautological bundles, and $\widehat{\mathscr O}_{\mathrm{vir}}$ is the symmetrized virtual structure sheaf. The subscript $\mathrm{loc}$ means that equivariant characters needed for localization have been inverted; $\mathbf d$ runs over effective multidegrees and $z^{\mathbf d}=\prod_i z_i^{d_i}$.
\begin{Def}
The element
\begin{equation}
V^{(\tau)}(z)=\sum\limits_{\textbf{d}=0}^{\infty} z^{\textbf{d}} {\rm{ev}}_{p_2, *}\Big(\textsf{QM}^{\textbf{d}}_{{\rm{nonsing}} \, p_2},\widehat{{\mathscr{O}}}_{{\rm{vir}}} \tau (\left.\mathscr{V}_i\right|_{p_1}) \Big) \in  K_{\mathsf{T}}(X)_{loc}[[z]]
\label{eq:vertexQKgen}
\end{equation}
is called bare vertex with descendent $\tau\in K_T(X)$.
\end{Def}
The vertex function is an element of localized quantum K-theory of $X$. Using the localization theorem one can compute the above vertex functions as equivariant integrals over the moduli space of quasimaps \cites{Pushkar:2016qvw,Koroteev:2021tl}. The ring structure of the K-theory arises in the saddle point limit of the vertex function in the equivariant parameter $q$ corresponding to the $\mathbb{C}^\times$-action on the base curve.
In the next formula $\mathbf p$ labels a torus fixed point of $T^*\mathbb F l_n$, not a momentum. The $x_{i,j}$ are the weights of the rank-$i$ tautological bundle at this fixed point, $d_i=\sum_{j=1}^i d_{i,j}$, and $d_{n,k}=0$. We use $(u,q)_d=\prod_{b=0}^{d-1}(1-uq^b)$ for $d\geq0$ and $(u,q)_{-d}=1/(uq^{-d},q)_d$. The cone $C$ consists of the nonnegative degree arrays contributing to the fixed-point quasimap localization formula.
\begin{Prop}
For $X$ being a cotangent bundle to complete flag variety of type A the vertex function coefficient for the identity class $\tau=1$ reads
\begin{equation}
\begin{aligned}
V^{(1)}_{\textbf{p}}(z)
={}&\sum_{(d_{i,j})\in C}\prod_{i=1}^{n-1}\Bigg[
 \left(\frac{q}{\hbar}\frac{\xi_i}{\xi_{i+1}}\right)^{d_i}
 \prod_{j,k=1}^i
 \frac{(q x_{i,j}/x_{i,k},q)_{d_{i,j}-d_{i,k}}}
      {(\hbar x_{i,j}/x_{i,k},q)_{d_{i,j}-d_{i,k}}}\\[-2pt]
 &\hspace{35mm}{}\times
 \prod_{j=1}^i\prod_{k=1}^{i+1}
 \frac{(\hbar x_{i+1,k}/x_{i,j},q)_{d_{i,j}-d_{i+1,k}}}
      {(q x_{i+1,k}/x_{i,j},q)_{d_{i,j}-d_{i+1,k}}}
 \Bigg].
\end{aligned}
\label{eq:V1pdef}
\end{equation}
where we define $x_{n,k}=a_k$ -- equivariant parameters for the final node of the quiver, and $z_i=\frac{\xi_{i}}{\xi_{i+1}}$ are the quantum (K\"ahler) parameters.
\end{Prop}
In \cites{Koroteev:2018azn,Koroteev:2021aa} it was shown that normalized quantum K-theoretic vertex functions for cotangent bundles to flag varieties of type A are the eigenfunctions of the tRS Hamiltonians. 

We use the characteristic-polynomial normalization of the Macdonald operators. For a multiplicative coordinate tuple $u$, let $p_i$ act by $u_i\mapsto qu_i$. With interaction parameter $t$ the operators are
\begin{equation*}
 t^{\binom r2}\sum_{|I|=r}\prod_{i\in I,\,j\notin I}
 \frac{tu_i-u_j}{u_i-u_j}\prod_{i\in I}p_i,\qquad 1\le r\le n.
\end{equation*}
Taking $(u,t)=(\xi,\hbar)$ gives $\mathcal H_r^{tRS}$ in the magnetic frame. For the electric operator $\mathcal H_r^{!\,tRS}$ one takes $(u,t)=(a,q/\hbar)$ in the vertex convention of \cite{Koroteev:2018azn}. Its classical interaction tends to $\hbar^{-1}$. A common rescaling of the spectral tuple fixes the central character. In particular the centered convention, with overall factor $t^{-r(n-r)/2}$ in the displayed sum, passes to ours by multiplying every spectral variable by $t^{(n-1)/2}$.

\begin{Thm}[\cite{Koroteev:2021aa}]\label{Th:tRSeqproof}
After the fixed-point normalization by a multiplicative prefactor, the vertex functions $\mathsf V^{(1)}_{\mathbf p}$ satisfy
\begin{equation}
\mathcal H_r^{tRS}(\xi) \mathsf{V}^{(1)}_{\textbf{p}} = e_r (\mathbf{a}) \mathsf{V}^{(1)}_{\textbf{p}}\,, \qquad r=1,\dots, n\,,
\label{eq:tRSEigenz}
\end{equation}
where $e_r$ is the elementary symmetric polynomial of degree $r$.
\end{Thm}

Similarly, in the electric frame
\begin{Thm}[\cite{Koroteev:2018azn}]\label{Th:tRSeqproofel}
With the normalization for shifts in the equivariant variables and the corresponding choice of the common spectral character, the vertex functions satisfy
\begin{equation}
\mathcal H_r^{!\,tRS}(a) \mathsf{V}^{(1)}_{\textbf{p}} = e_r (\xi) \mathsf{V}^{(1)}_{\textbf{p}}\,, \qquad r=1,\dots, n\,,
\label{eq:tRSEigena}
\end{equation}
where the successive ratios of the $\xi_i$ are the K\"ahler parameters.
\end{Thm}
The arguments use kernel identities of the kind studied in \cites{2013arXiv1305.4759H,2012arXiv1206.3787H,2012JMP....53l3512H}.

Thus the normalized vertex functions solve the tRS spectral problem. At $q\to1$ the classical characteristic coefficients give relations in the quantum K-theory ring of the quiver variety. At a root of unity the shift has finite order. The quantum K-theory analysis with A. Smirnov \cite{KSRoots} and the monodromy eigenvectors of Tarasov and Varchenko \cite{TVRoots} lead to a related statement for complete flags. The Newton--Adams combinations of the tRS Hamiltonians reproduce the Frobenius action in every rank and at every root order \cite{BetheCycles}.
\begin{Thm}[\cite{Koroteev:2021tl}]
For a type-$A_r$ quiver variety on the nondegenerate Bethe locus, the eigenvalues of quantum multiplication by the quantum tautological class $\widehat\tau(z)$ are the characters $\tau(s_I)$ evaluated at solutions of the type-$A_r$ XXZ Bethe equations.
\end{Thm}

The generating functions of quantum exterior powers of the tautological bundles give Baxter operators. The dimension vector of the quiver variety selects the corresponding weight space in the spin-chain representation of the quantum affine algebra of type $A_r$ \cites{Pushkar:2016qvw,Koroteev:2021tl}.
\begin{Thm}[\cite{Koroteev:2021tl}]
On the regular chart, the localized presentation for the quantum equivariant K-theory of $T^*\mathbb F l_n$, after adjoining a common multiplicative coordinate to the K\"ahler ratios, reads
\begin{equation}
QK_T(T^*\mathbb F l_n)_{\rm loc}\simeq
\frac{\mathbb C[\zeta_i^{\pm1},a_i^{\pm1},p_i^{\pm1},\hbar^{\pm1}\mid 1\le i\le n]_{\rm loc}}
     {(\mathcal H_r^{tRS}(\zeta,p,\hbar)-e_r(a))_{r=1}^{n}}\,.
\label{eq:Kthflag}
\end{equation}
\end{Thm}

Here the $p_j$ commute. Geometrically they are parameter-normalized quantum products of the determinant of the $j$th tautological bundle and the dual determinant at the preceding node \cite{Koroteev:2021tl}. The formal quantum K-ring uses completion in the effective K\"ahler degrees.

\subsection{Truncations to Macdonald Polynomials}
Macdonald polynomials $P_\lambda(\xi;q,\hbar)$ can be obtained from infinite series \eqref{eq:V1pdef} by specifying equivariant parameters $a_1,\dots,a_n$ on a lattice spanned by $q$ and $\hbar$.
\begin{Prop}[\cite{Koroteev:2021aa}]
\label{Prop:Truncation}
Let $\lambda_1\ge\dots\ge\lambda_n\ge0$ be a partition, with zero parts appended as needed. In the magnetic normalization, specialize the unordered spectral parameters to
\begin{equation}
\{a_1,\ldots,a_n\}=\{q^{\lambda_i}\hbar^{n-i}\mid 1\le i\le n\}.
\label{eq:truncationa}
\end{equation}
For a compatible fixed point $\mathbf q$ and chamber the bare series terminates. Restoring the leading monomial and choosing its coefficient to be one gives
\begin{equation}
\mathsf{V}^{(1)}_{\textbf{q}} = P_\lambda(\xi;q,\hbar)\,.
\end{equation}
\end{Prop}

\subsection{The DELL System}
Together with Shakirov, I proposed quantum Hamiltonians for the double elliptic (DELL) system in oval 1 \cite{Koroteev:2019gqi}. The dependence on both momenta and positions is elliptic. Grekov and Zotov expressed the generating function by a characteristic determinant and constructed its classical matrix with a Manakov triple \cite{Grekov:2020jpk}. Their large-$N$ construction is developed in \cite{Grekov:2021zqq}.

For $N$ particles the center reduction leaves $N-1$ degrees of freedom. The proposed Hamiltonians are
\begin{equation}
\widehat{H}_a = \mathcal{O}_0^{-1} \mathcal{O}_a\,, \qquad a = 1, \ldots, N-1\,,
\label{eq:DELLHAMs}
\end{equation}
where operators $\mathcal{O}_0,\mathcal{O}_1,\dots,\mathcal{O}_{N-1}$ are Fourier modes of the following current
\begin{equation}
\mathcal{O}(z) \ = \ \sum\limits_{n \in {\mathbb Z}} \ \mathcal{O}_n \ z^n \ = \ \sum\limits_{n_1, \ldots, n_N = -\infty}^{\infty} \ (-z)^{\sum n_i} \ w^{\sum \frac{n_i(n_i - 1)}{2}} \ \prod\limits_{i < j} \theta_p\big( \hbar^{n_i - n_j} x_i /  x_j \big) \ p_1^{n_1} \ldots  p_N^{n_N}\,.
\label{eq:Ocurrent}
\end{equation}
Here $p$ and $w$ are the elliptic parameters for the coordinates and momenta, respectively. We use the multiplicative theta factor of \cite{Koroteev:2019gqi} and place all coefficients to the left of the shifts. Explicitly $\widehat p_i f(x)=f(x_1,\ldots,qx_i,\ldots,x_N)$, so $\widehat p_i\widehat x_j=q^{\delta_{ij}}\widehat x_j\widehat p_i$. We interpret the series formally and invert $\mathcal O_0$ in the formal localization.
Pairwise commutativity of the full quantum family is Conjecture~1.1 of \cite{Koroteev:2019gqi}, where low-rank checks are given.

A second conjecture relates the joint eigenfunction to the elliptic genus of the affine Laumon spaces. In the defect partition function one sums over the degrees and takes the normalized Nekrasov--Shatashvili limit as in \cite{Koroteev:2019gqi}.
\begin{Con}[\cite{Koroteev:2019gqi}]\label{th:Conj2}
Let $\mathcal{Z}^{DELL}_{\text{inst}}(p,x_1,\dots,x_N)$ denote the resulting normalized generating function, with its equivariant parameters understood. Then there exists a function $\lambda(z,\textbf{a},w,p)$ such that
\begin{equation}
\mathcal{O}(z)\mathcal{Z}^{DELL}_{\text{inst}}(p,x_1,\dots, x_N) = \lambda(z, \textbf{a},w,p) \ \mathcal{O}_{0} \mathcal{Z}^{DELL}_{\text{inst}}(p,x_1,\dots, x_N)\,.
\label{eq:Oeignproblem}
\end{equation}
Taking the coefficient of $z^n$ and applying $\mathcal O_0^{-1}$ gives the eigenvalue equation for $1\le n<N$, where $\lambda_n$ is the corresponding coefficient of $\lambda$:
\begin{equation}
\widehat H_n \mathcal{Z}^{DELL}_{\text{inst}}(p,x_1,\dots, x_N) = \lambda_n(\textbf{a},w,p) \mathcal{Z}^{DELL}_{\text{inst}}(p,x_1,\dots, x_N)\,.
\end{equation} 
\end{Con}
The computations in \cite{Koroteev:2019gqi} give evidence for the eigenfunction conjecture. At $N=2$ there is one reduced Hamiltonian, so quantum commutativity imposes its first nontrivial test at $N=3$.

Braden, Marshakov, Mironov and Morozov constructed the classical two-particle double elliptic Hamiltonian by a duality argument \cite{BMMM}. Subsequent work relates the spectral curve to six-dimensional gauge theory compactified on a torus and to the Seiberg--Witten curve in M-theory \cites{Hollowood:2003cv,Braden:2001yc,Braden:2003gv}. The two elliptic parameters come from the gauge theory and the compactification torus.

The self-dual constructions and their relation to the KS operators are discussed in \cite{Mironov:2023ab}. In particular the two upper degenerations have different meanings. Sending $w$ to zero gives ordinary eRS (oval 3), whereas degenerating the coordinate elliptic parameter gives a model elliptic in the momenta (oval 2). The scaling limits with varying coupling are described in \cite{Gorsky:2022aa}. Awata et al. propose explicit series for eigenfunctions of the KS operators and compare them with instanton calculations \cite{AwataDell}.

\subsection{The Elliptic RS Model}
The eRS model from oval \# 3 can be understood as either elliptic deformation of the tRS model or as the limit $w\to 0$ of the DELL system above. Its Hamiltonians read
\begin{equation}\label{eq:eRSEigenproblemTh}
\mathcal{H}^{eRS}_r =\sum_{\substack{\mathcal{I}\subset\{1,\dots,N\} \\ |\mathcal{I}|=r}}\prod_{\substack{i\in\mathcal{I} \\ j\notin\mathcal{I}}}\frac{\theta_1(\hbar x_i/x_j|\mathfrak{p})}{\theta_1(x_i/x_j|\mathfrak{p})}\prod\limits_{i\in\mathcal{I}}p_i \,,
\end{equation}
Here $0<|\mathfrak p|<1$ is the coordinate elliptic parameter, $\theta_1$ is the odd multiplicative theta function, and $p_i$ shifts $x_i$ by $q$.

There are two questions about the wave function: its construction as an eigenfunction and its identification with a particular enumerative series. Langmann, Noumi and Shiraishi construct convergent eigenfunctions perturbatively in appropriate domains \cite{EllipticRSWave}. The geometric identification can be formulated as follows.
\begin{Con}[\cites{Bullimore:2015fr,ShiraishiAffine}]\label{con:affineRS}
A joint wave function for the elliptic Ruijsenaars--Schneider operators,
\begin{equation}\label{eq:eRSEigenproblem}
\mathcal{H}^{e\text{RS}}_k \mathcal{Z}^{{\text{RS}}}(\textbf{a},\textbf{x}) = \lambda_k (\textbf{a})\mathcal{Z}^{{\text{RS}}}(\textbf{a},\textbf{x})\,,\qquad k=1,\dots,N-1\,.
\end{equation}
is obtained, with the prescribed prefactor and stationary limit, from the generating series of equivariant K-theoretic Euler characteristics of the affine Laumon spaces $\mathcal L^{\rm aff}_{\mathbf d}$. We write the series schematically as
\begin{equation}
\sum_{\mathbf d\ge0}\mathfrak q^{\mathbf d}\,
 \chi_T(\mathcal L^{\mathrm{aff}}_{\mathbf d},\,\textstyle\sum_j(-y)^j\Omega^j),
\label{eq:equivaraintKthLaumon}
\end{equation}
where $y$ records the exterior powers in the equivariant $\chi_y$-genus, and the components of $\mathfrak q$ weight the degrees. The stationary prescription and parameter conventions are given in \cite{ShiraishiAffine}.
\end{Con}

For comparison, Awata et al. construct a nonstationary difference operator in arbitrary rank and conjecture its affine Laumon solution \cite{LaumonIII}. Di Francesco and Kedem prove a nonstationary difference equation in the affine relativistic Toda specialization \cite{LaumonToda}.
\begin{Thm}[\cite{Gorsky:2022aa}]\label{RSSpec1}
Let $\mathcal Z^{\mathrm{RS}}$ be a joint eigenfunction of \eqref{eq:eRSEigenproblemTh} regular at the evaluation points below and nonzero at $\hbar^\rho$. Its eigenvalues satisfy
\begin{equation}
\lambda_k(\textbf{a})=\prod\limits_{n = 0}^{k-1} \frac{\theta_1(\hbar^{N-n})}{\theta_1(\hbar^{n+1})} \cdot\frac{\mathcal{Z}^{{\text{RS}}}(\textbf{a},\hbar^{\rho}q^{\omega_k})}{\mathcal{Z}^{\text{RS}}(\textbf{a},\hbar^{\rho})}\,,\qquad k=1,\dots, N-1\,,
\end{equation}
Here $\rho_i=(N+1)/2-i$ and $q^{\omega_k}$ acts by multiplying the first $k$ coordinates by $q$. At $x=\hbar^\rho$ only the subset $\{1,\ldots,k\}$ survives in the $k$th Hamiltonian, giving the formula by evaluation. Applying it to the affine Laumon series uses the identification in Conjecture~\ref{con:affineRS}.
\end{Thm}

\subsection{The Calogero-Moser Systems}
For ordinary Laumon spaces, Negu\c{t} proved Braverman's conjecture relating equivariant integrals of the Chern polynomial of the tangent bundle to the trigonometric Calogero--Sutherland eigenfunction \cites{Braverman:aa,Negut_2009}. His affine Laumon construction treats the elliptic Calogero--Moser problem with a full surface operator \cite{NegutAffine}.

There is a complementary gauge-theoretic description through $qq$-characters \cites{Nekrasov:2015wsu,Nekrasov:2017gzb}. Grekov and Nekrasov obtain bilinear identities for the elliptic Calogero--Moser system from crossed and folded instantons \cite{GrekovNekrasov}. Grekov derives the classical Lax matrix and its eigenvector from instanton moduli and describes K-theoretic and elliptic counterparts \cite{GrekovInstantons}.

\subsection{Spaces of Opers for Elliptic Many-Body Systems}
On the affine Laumon diagonal the spectral curve and its sheaf supply geometric data beyond the top Hamiltonian level. For ordinary eRS (oval 3) and eCM (oval 6), one evaluates a Baker--Akhiezer section on a marked spectral cover. Shifted evaluations give the difference Wronskians in the former case, and derivatives give the differential Wronskians in the latter. Specifying the auxiliary spectral line and the evaluation divisor fixes the flag. On the reflected side, the rank-two elliptic Wronski map and the dynamical elliptic Bethe algebra were studied by Thompson and Varchenko \cite{TVell}. Identifying the monodromy data with a prescribed dual eCM level in oval 4 is an additional classical problem.

At the DELL vertex the base of the proposed construction is an abelian surface. On its cyclic cover the translates of a theta section define a marked flag on the regular locus. The determinant descends as a theta section of polarization type $(1,N)$ on the cyclic quotient. Dividing by the fixed intertwiner determinant gives the classical characteristic current of \cite{Grekov:2020jpk}. In the split construction, fixing the ordered flat characteristics and the starting row gives a finite compatibility test for prescribed theta minors in arbitrary rank. The test reconstructs a common set of columns when it succeeds. A fixed quantum lattice Hamiltonian realizing these Bethe data remains to be constructed. Even on a two-particle trigonometric level the cyclic and Baker--Akhiezer first flags can have different degrees, so matching the classical characteristic polynomial alone does not identify their markings.

As far as space of opers and dual integrable systems are concerned, another interesting direction is to extend the setup to superalgebras. Some progress in this direction has been made in \cites{Zeitlin:2013aa,Huang:2018aa,Zeitlin:2023ab}.

\subsection{DAHAs}\label{Sec:DAHA}
For type $A$, one description of DAHA begins with the braid group of configurations of particles on a torus, with its central parameter $q$. The two cycles give the $X$ and $Y$ generators. The braid generators $T_i$ satisfy quadratic Hecke relations, for example $(T_i-\hbar^{1/2})(T_i+\hbar^{-1/2})=0$ in the symmetric normalization. This gives the algebra at oval 5 \cites{Cherednik-book,GKNS,Matsuo:2023aa}.

The punctured torus also enters the gauge-theoretic construction of the character variety for the $\mathcal N=2^*$ theory \cite{GKNS}.

If $e$ is the symmetrizing idempotent of the finite Hecke algebra, the spherical DAHA is the corner $e\,\mathrm{DAHA}\,e$. In type $A$ its classical limit is related to the multiplicative Calogero--Moser space. In rank one it quantizes the character variety of flat $SL(2,\mathbb C)$ connections on the once-punctured torus with fixed puncture trace \cites{Oblomkov:aa,oblomkov2004double,GKNS}.

In the polynomial representation the commuting spherical Macdonald operators realize the quantum tRS Hamiltonians. Cherednik's difference-elliptic and differential elliptic operators are relevant to ovals 3 and 6 \cites{elldifdaha,Cherednik1995}. Chalykh constructs quantum Lax pairs for elliptic Calogero--Moser and generalized Ruijsenaars systems using Dunkl and Cherednik operators \cite{ChalykhLax}. Rains constructs elliptic double affine Hecke algebras as sheaf algebras \cite{RainsEll}. Identifying such an algebra with a particular oval requires a choice of parameters and representation.

In the polynomial representation the spherical DAHA acts on symmetric Laurent polynomials in $n$ variables over $\mathbb Q(q,\hbar)$. At the specialization \eqref{eq:truncationa}, the normalized vertex functions give Macdonald polynomials with the eigenvalues \eqref{eq:tRSEigenz}.

The symmetric Macdonald polynomials are indexed by partitions of length at most $n$. For generic parameters they form a simultaneous eigenbasis for the commuting Macdonald operators. The action of symmetric multiplication on this basis is described by the Pieri rules \cite{Cherednik-book}.

The degenerations in ovals 7--9 lead to trigonometric and rational Cherednik algebras. In rank one, the brane quantization of the punctured-torus character variety provides a concrete comparison with the spherical DAHA. In \cite{GKNS} compactly supported Lagrangian $A$-branes are matched with finite-dimensional modules in explicit families, including their shortening conditions and examples of extensions. Claim~1.1 of that work proposes a derived equivalence on the corresponding subcategories. Extending the categorical comparison to higher rank and constructing an algebra whose commuting operators recover the full DELL family remain open questions for the diamond.


\bibliography{tokyo-review}
\end{document}

%% file: oval-key.tex
\begingroup
\definecolor{ovalink}{HTML}{27333A}
\definecolor{ovalphase}{HTML}{17549B}
\definecolor{ovalframe}{HTML}{785B12}
\definecolor{ovalfill}{HTML}{FCF8E7}
\begin{tikzpicture}[x=1pt,y=1pt,every node/.style={inner sep=0pt,text=ovalink}]
  \path[draw=ovalframe,fill=ovalfill,line width=.65pt,rounded corners=23pt]
    (-89,-46.5) rectangle (89,46.5);
  \node[font=\fontsize{10}{12}\selectfont] at (0,33)
    {\textcolor{magenta}{\textbf{\#}}\enspace\textbf{Model}\enspace
     \textcolor{ovalphase}{momentum $\times$ position}};
  \node[font=\fontsize{9}{11}\selectfont] at (0,16) {Type of oper};
  \node[font=\fontsize{9}{11}\selectfont] at (0,2) {Magnetic model $\cdot$ Wronskian relation};
  \draw[ovalframe,opacity=.55,line width=.32pt] (-61.5,-7.5) -- (61.5,-7.5);
  \node[font=\fontsize{9}{11}\selectfont] at (0,-18) {Algebra for $n$ particles};
  \node[font=\fontsize{9}{11}\selectfont] at (0,-33) {Algebra as $n\to\infty$};
\end{tikzpicture}
\endgroup